\documentclass[reqno, 12pt]{amsart}
\usepackage{amsfonts}
\newcommand{\ol}{\overline}
\newcommand{\tG}{\Gamma}
\renewcommand{\phi}{\varphi}

\newcommand{\Ad} {\mathop{\rm Ad}}
\newcommand{\ad} {\mathop{\rm ad}}
\newcommand{\End}{\mathop{\rm End}}
\newcommand{\Id} {\mathop{\rm Id}}
\newcommand{\SSp} {\mathop{\rm Sp}}
\newcommand{\Hom} {\mathop{\rm Hom}}
\newcommand{\Vect} {\mathop{\rm Vect}}
\newcommand{\ssq}{\mathop{\rm S^2}}
\newcommand{\iso}{\cong}                    

\newcommand{\Lsl}{{\mathfrak{sl}}}
\newcommand{\Lg} {{\mathfrak g }}
\newcommand{\Lh} {{\mathfrak h }}
\newcommand{\Ly} {{\mathfrak y }}
\newcommand{\Lz} {{\mathfrak z }}
\newcommand{\Lm} {{\mathfrak m }}

\newcommand{\bbC}{\mathbb C}

\newcommand{\wrapX}{\widetilde X} 

\newcommand{\CM}{\mathcal{M}}
\newcommand{\CN}{\mathcal{N}}
\newcommand{\CE}{\mathcal{E}}
\newcommand{\CH}{\mathcal{H}}
\newcommand{\CF}{\mathcal{F}}

\newtheorem*{Theor*}{Theorem}
\newtheorem{Lemma}{Lemma}
\newtheorem*{Lemma*}{Lemma}
\newtheorem*{MainLemma*}{Basic Lemma}
\newtheorem*{Remark*}{Remark}

\title[Invariant linear connections on homogeneous symplectic varieties]
{Invariant linear connections\\ on homogeneous symplectic varieties}

\author[S.V.~Pikulin, E.A.~Tevelev]{S.V.~Pikulin, E.A.~Tevelev$^*$}

\begin{document}

\begin{abstract}
We find all homogeneous symplectic varieties
of connected reductive algebraic groups that
admit an invariant linear connection.
\end{abstract}

\maketitle

\thispagestyle{empty}
\setcounter{footnote}{1}
\footnotetext{
The research was supported by the grant
INTAS-OPEN-97-1570 of the INTAS foundation.}

\section*{Introduction}
Let $G$~be a connected Lie group.
The following problem was stated many years ago:\break
Is it possible to describe all homogeneous $G$-spaces $X$ that admit an
invariant linear connection 
(see~\cite{MMO}, \cite{Wang}, \cite{Nomizu})?
Not much is known in general and the main attention
was focused on the description of invariant linear connections
on some particularly nice homogeneous spaces, namely, on symmetric
spaces or more generally on reductive homogeneous spaces 
(see~\cite{Koba}, \cite{Kur}). The aim of this paper
is to give the complete solution of this problem in the following
situation: $G$ is a connected reductive algebraic group over $\mathbb C$,
$X$ is a symplectic $G$-variety, that is, $X$ is either
an adjoint orbit of $G$
or its covering.

We shall explain in \S1 that it is sufficient to consider
only simple simply-connected groups $G$ and only their nilpotent
orbits $X$ (and their coverings). In this set-up our answer
is given by the following theorem:

\begin{Theor*} Suppose that $G$ is a simple simply-connected group,
$X$ is its nilpotent orbit. The following assertions are equivalent:
\begin{itemize}
\item[(a)] there exists an invariant linear connection on $X${\em;}
\item[(b)] there exists an invariant linear connection on some covering 
$\wrapX$ of $X${\em;}
\item[(c)] $G\iso\SSp_n$ and $X=(\Ad G)e$, where $e$ is a
highest root vector.
\end{itemize}
\end{Theor*}

The paper is organized as follows. In \S1 we shall remind
some basics on invariant linear connections and prove
some auxiliary lemmas. In \S3 we shall prove the Theorem.
First in \S2 we shall prove the Basic Lemma, which says that
if there exists an invariant linear connection on the nilpotent orbit~$X$
then $X$ is locally $G$-equivariantly isomorphic to an open
orbit in some $G$-module. The Basic Lemma is proved a~priori,
using some $\Lsl_2$-trick. Then we shall prove the Theorem
using the case-by-case considerations. It is possible (and easy),
since the list of all $G$-modules having an open orbit is well-known
(and short).

The authors are grateful to their supervisor E.B.~Vinberg
for useful discussions and the simplification of some proofs.

\section*{\S1. Basics on invariant linear connections}

Let $G$ be a connected algebraic group over $\mathbb C$
with Lie algebra $\Lg$. Let $X$~be a homogeneous $G$-variety.
Denote by $H=G_o\subset G$~the stabilizer of some point $o\in X$, 
let $\Lh$ be a Lie algebra of $H$. 
The tangent space $T_o X$ is canonically identified with $\Lg/\Lh$.
For any $\xi\in\Lg$ we denote the vector $\xi+\Lh\in \Lg/\Lh$ by $\ol\xi$.
If $X:\Lg\to\Lg$~is a linear operator such that
$X\Lh\subset\Lh$ then $\ol X$ is the induced operator on $\Lg/\Lh$,
$\ol X(\ol\xi)=\ol{X\xi}$.
Any $\xi\in\Lg$ induces the vector field $\xi^*$ on $X$.
The corresponding map $\Lg\to\Vect(X)$ is a homomorphism of Lie algebras.
Since $X$ is homogeneous, the values of vector fields $\xi^*$ 
at any point of $X$ span the tangent space at this point.
In particular, we have $\xi^*(o)=\ol\xi$.

It is well-known (see e.g.~\cite{MMO}, Theorem~2)
that there exists a 1-1 correspondence between
the set of invariant linear connections on $X$ 
and the set of linear maps
$\Gamma:\Lg\to\End {\Lg/\Lh}$ that satisfy the following conditions:
$$\hbox{\rm I)}\quad \Gamma_\alpha=\ol{\ad\alpha}\ \hbox{\rm for any}\ \alpha\in\Lh;$$
$$\hbox{\rm II)}\quad \Gamma_{\Ad(x)\xi}=\ol{\Ad(x)}\Gamma_\xi\ol{\Ad(x)}^{-1}
\ \hbox{\rm for any}\ x\in H,\ \xi\in\Lg.$$
If $H$ is connected then II) is equivalent to
$$\hbox{\rm II')}\quad \Gamma_{[\alpha,\xi]}=[\Gamma_\alpha,\Gamma_\xi]
\ \hbox{\rm for any}\ \alpha\in\Lh,\ \xi\in\Lg.$$

The covariant derivative $\nabla_{\xi^*}\eta^*(o)$ of 
the vector field $\eta^*$ in the direction of the vector field $\xi^*$ 
at the point $o$ is given by
$$
\nabla_{\xi^*}\eta^*(o)=\Gamma_\xi\ol{\mathstrut\eta}+\ol{[\xi,\eta]}.
$$
The curvature and the torsion of the connection $\nabla$ at the point $o$ 
are given by 
\begin{align}
\rho(\ol{\mathstrut\xi},\ol{\mathstrut\eta})=
\Gamma_{[\xi,\eta]}-[\Gamma_\xi,\Gamma_\eta],
\notag \\
\sigma(\ol{\mathstrut\xi},\ol{\mathstrut\eta})=
\Gamma_\xi\ol{\mathstrut\eta}-\Gamma_\eta\ol{\mathstrut\xi}-
\ol{[\xi,\eta]}.
\label{sigma}
\end{align}
A connection is called locally flat if it has zero curvature and torsion.

It is important that we can average connections 
in the following sense:

\begin{Lemma}
Suppose that $K\subset G$ is a reductive subgroup that normalizes $H$. 
Suppose that there exists an invariant linear
connection on $X$. Then there exists an invariant linear connection
such that the corresponding map $\Gamma$ satisfies the following 
additional condition:
$$\Gamma_{\Ad(x)\xi}=\ol{\Ad(x)}\Gamma_\xi\ol{\Ad(x)}^{-1}
\ \hbox{\rm for any}\  x\in K,\ \xi\in\Lg.$$
\end{Lemma}

\begin{proof}
We note that conditions I) and II) define the {\em affine} subspace $M$
in the linear space
$\CM=\Hom(\Lg,\End {\Lg/\Lh})$.
Since there exists an invariant linear
connection on $X$, $M$ is non-empty.
Clearly $\CM$ is a $K$-module and we need to check that $M$
contains a $K$-fixed element.
But this is evident since $M$ is $K$-invariant and $K$ is reductive.
Namely, let $\CN\subset\CM$ be a linear subspace
corresponding to the affine subspace~$M$.
Then $\CN$ is a $K$-submodule and there exists a complementary
$K$-submodule $\CN^{\perp}$.
The intersection $\CN^{\perp}\cap \CM$ is fixed by $K$
and represents the desired connection.
\end{proof}

From now on we suppose that $X$ is an adjoint orbit, $X=\Ad(G)o$.
At this point $G$~still can be non-reductive. 
So $H=\{g\in G\,|\,\Ad(g)o=o\}$~is the centralizer of $o$ in $G$ and
$\Lh=\Lz(o)=\{\gamma\in\Lg\,|\,[\gamma,o]=0\}$.
First we want to perform a reduction to nilpotent orbits
promised in the Introduction.
Take the Jordan decomposition $o=o_s+o_n$, where
$o_s$~is semisimple, $o_n$~is nilpotent, and $[o_s,o_n]=0$.
Let $Z(o_s)\subset G$ be the centralizer of $o_s$.
It~is~well-known that $Z(o_s)$~is a connected subgroup of $G$. 
If $G$ is reductive then 
$Z(o_s)$ is reductive as well. Its Lie algebra is
$\Lz(o_s)=\{\xi\in\Lg\,|\,[\xi,o_s]=0\}$. We have
$o_n\in\Lz(o_s)$, $\Lh=\Lz(o_s)\bigcap\Lz(o_n)$.

\begin{Lemma}
There exists an invariant linear
connection on $X$ (resp.~on its covering) iff there exists a
$Z(o_s)$-invariant
linear connection on $(\Ad Z(o_s))o_n$ (resp.~on its covering).
\end{Lemma}

\begin{proof}
We will prove only the assertion about orbits, leaving
the assertion about coverings to the interested reader.
Though $Z(o_s)$ is not necessarily reductive, 
the homogeneous space $G/Z(o_s)$ is always reductive,
that is, $\Lg=\Lz(o_s)\oplus\Lm$, 
where $\Ad(Z(o_s))\Lm\subset\Lm$.
Suppose that $\Gamma^{(s)}:\Lz(o_s)\to\End\Lz(o_s)/\Lh$~is a linear map
satisfying I) and~II) for $x\in Z(o_s)$, $\xi\in\Lz(o_s)$.
Then we can extend $\Gamma^{(s)}$ to the linear map
$\Gamma:\Lg\to\End\Lg/\Lh$ satisfying I) and~II) for $x\in G$, $\xi\in\Lg$.
Namely, we take $\Gamma_\Lm=0$ and $\Gamma_{\Lz(o_s)}|_\Lm=\ol{\ad\Lz(o_s)}$.
Conversely, if we have an appropriate linear map
$\Gamma:\Lg\to\End\Lg/\Lh$ then we can define a linear map
$\Gamma^{(s)}:\Lz(o_s)\to\End\Lz(o_s)/\Lh$ by taking a projection
of $\End\Lg/\Lh$ on the subspace of operators 
preserving $\Lz(o_s)$ and $\Lm$ and by restricting $\Gamma$ on $\Lz(o_s)$.
\end{proof}

So from now on  $o=o_n$~is a nilpotent element.
Certainly it will be denoted by $e$.
In the sequel we suppose that $G$ is reductive.
Clearly we can assume also that $G$ is semisimple and simply-connected,
since this assumption does not change the set of nilpotent orbits.
It is not hard to see that we can actually suppose that $G$ 
is simple by taking each simple component.

By the Jacobson--Morozov theorem
the nilpotent element $e\in\Lg$ can be included into the
$\Lsl_2$-triple $(f,h,e)\subset\Lg$.
Let $T$ be a 1-PS such that its Lie algebra is spanned by $h$.
Clearly $T$ normalizes $H=G_e$ and hence its identity component $H^0$.
We have

\begin{Lemma}
Suppose that there exists an invariant linear
connection on the universal covering $\wrapX=G/H^0$ of $X$.
Then there exists an invariant linear connection on $\wrapX$
such that the corresponding map~$\Gamma$ satisfies the following 
additional condition:
\begin{align}
\Gamma_{[h,\xi]}=[\ol{\ad h},\Gamma_\xi]\ \hbox{\rm for any}\ 
\xi\in\Lg.
\label{adh}
\end{align}
\end{Lemma}

\begin{proof}
This is an easy application of Lemma~1.
\end{proof}

\section*{\S2. Invariant linear connections on nilpotent orbits }

\begin{MainLemma*} 
Suppose that there exists an invariant linear
connection on $\wrapX$.
Then there exists an invariant linear locally flat connection
on $\wrapX$ such that 
\begin{align}
\Gamma_\xi\ol h=\ol\xi\ \hbox{\rm for any}\  \xi\in\Lg.
\label{hrep}
\end{align}
\end{MainLemma*}

It follows from~\cite[\S5]{MMO} that the invariant linear connection
$\nabla$ given by the linear map $\Gamma:\Lg\to\End {\Lg/\Lh}$
is locally flat if and only if the map $\Gamma$ 
is a representation of $\Lg$ in ${\Lg/\Lh}$, which is
transitive. That is, there exists a vector $v\in \Lg/\Lh$
such that $\Gamma_\Lg v=\Lg/\Lh$. So let us show that $\Gamma$
given by Lemma~3 is a transitive representation of $\Lg$.

\begin{Lemma} 
There exists an $\Lsl_2$-triple $\langle\CF,\CH,\CE\rangle\subset\End {\Lg/\Lh}$,
such that $\CE=\ol{\ad e}$\break and $\CH=\ol{\ad h}+\Id_{\Lg/\Lh}$.
\end{Lemma}
\begin{proof}
Let $Y\subset\Lg$~be an irreducible $\Lsl_2$-submodule
in $\Lg$ with the highest weight $\mu$ and the highest vector~$y$.
Let $\Ly\equiv Y\bigcap\Lh=\langle y\rangle $.
The operator $\ad e|_Y$ induces a structure of 
$\Lsl_2$-module on $\ol Y\iso Y/\Ly$.
The highest weight is $\mu-1$ and the highest vector is $\ol{[f,y]}$.
We notice that if we have $[h,\xi]=k\xi$ for $\xi\in Y$
then we have $\CH|_{\ol Y}\ol\xi=(k+1)\ol\xi$, therefore 
$\CH=\ol{\ad h}+\Id_{\Lg/\Lh}$.
\end{proof}

It follows from Lemma 4, from equation~\eqref{adh}, and from II') 
that for any $\xi\in\Lg$ we have
\begin{align}
\tG_{[e,\xi]}=[\CE,\tG_\xi],\quad
\tG_{[h,\xi]}=[\CH,\tG_\xi]\quad \hbox{\rm for any}\ \xi\in\Lg.
\label{H}
\end{align}

It is well-known that for any reductive group $S$
and any $S$-modules $U$ and $V$ the linear map $U\to V$
is $S$-equivariant if and only if it is equivariant 
with respect to the Borel subgroup of $S$.
It follows that the equations~\eqref{H} imply that 
we also have
\begin{align}
\tG_{[f,\xi]}=[\CF,\tG_\xi]\ \hbox{\rm for any}\ \xi\in\Lg.
\label{F}
\end{align}

Substituting $e$ and $h$ instead of $\xi$ in~\eqref{F} we get
$\CH=\tG_h$ and $\CF=\tG_f$.
Therefore,
\begin{align}
\tG_{[f,\xi]}=[\tG_f,\tG_\xi],\quad
\tG_{[h,\xi]}=[\tG_h,\tG_\xi],\quad
\tG_{[e,\xi]}=[\tG_e,\tG_\xi]\ \hbox{\rm for any}\ \xi\in\Lg.
\label{tG}
\end{align}

The subalgebras $\Lh\subset\Lg$ and 
$\langle f,h,e\rangle \subset\Lg$ together generate $\Lg$. 
Therefore, II') and~\eqref{tG} imply (by an easy induction) that
\begin{align}
\tG_{[\xi,\eta]}=[\tG_\xi,\tG_\eta]\ \hbox{\rm for any}\ \xi,\eta\in\Lg. 
\label{star}
\end{align}
So $\tG$ is a representation of $\Lg$ in ${\Lg/\Lh}$!

It remains to verify that $\tG$ is transitive.
We are going to check that $\tG_\Lg\ol f={\Lg/\Lh}$.
We have $\CF\ol f=0$ by the construction of $\CF=\tG_f$.
Moreover, vectors of the form $\CF^n\ol{[f,a]}$, 
where $a\in\Lh$, $n=0,1,\ldots$, linearly span ${\Lg/\Lh}$.
Let us prove that for any $n\ge0$ and $a\in\Lh$ 
the vector $\CF^n\ol{[f,a]}$ is equal to
$\tG_\xi\ol f$ for some $\xi\in\Lg$. 
Indeed, take $\xi=-(\ad f)^na$. Using~\eqref{star} we obtain
\begin{align*}
&\tG_\xi\ol f=-\tG_{[f,[f,\dots,[f,a]]\dots]}\ol f=
-[\tG_f,[\tG_f,\dots,[\tG_f,\tG_a]]\dots]\ol f=\\
&-\sum_{j=0}^n(-1)^j\binom n j\tG_f^{n-j}\tG_a\tG_f^j\ol f=
-\tG_f^n\tG_a\ol f=\tG_f^n\ol{[f,a]}=F^n\ol{[f,a]}.
\end{align*}
This proves that $\tG$ is transitive.

The formula~\eqref{hrep} follows from~\eqref{sigma} with $\eta=h$,
from $\sigma=0$ (since $\tG$ is transitive), and from the equality
$\tG_h=\CH=\ol{\ad h}+\Id_{\Lg/\Lh}$.
The Basic Lemma is proved.

\section*{\S2. Proof of the Theorem}

The implication $(a)\Rightarrow(b)$ is trivial, since
any invariant linear connection on $X$ 
can be lifted to an invariant linear connection on $\wrapX$.

Let us prove that $(c)\Rightarrow(a)$. Let $G\iso SSp_n$.
The adjoint representation of $G$ is isomorphic
to the symmetric square of the simplest representation of $G$ (cf.~\cite{VinO}).
This isomorphism takes the orbit of the highest root vector
to the orbit of non-zero ``perfect squares''
$\{u^2\in\ssq\bbC^n\,|\,u\in\bbC^{n}\}$.
This orbit is covered by the open orbit in the simplest representation
of $G$ in $\bbC^n$.
Therefore, $X$ is obtained from $\bbC^n\setminus\{0\}$ 
with a natural action of $G$ by
identifying opposite points.
Obviously this variety admits
an invariant linear connection.

Finally let us check that $(b)\Rightarrow(c)$.
It follows from $(b)$ that there exists an invariant linear connection
on the universal covering $\wrapX$.
Moreover, we can assume that this connection 
satisfies the assertions of the Basic Lemma.

It follows from~\cite[\S6]{MMO}
that all transitive representations of simple algebras $\Lg$ 
are given in the following table:
\begin{align*}
\begin{tabular}                                                {
| l | l          l             | c                            |}
\hline
    & $\Lg$      &             & $\tG$                        \\
\hline
  1 & $A_l$,     & $l\ge1$     & $kR(\pi_1),\ k=1,2,\dots,l$  \\
  2 & $A_{2l}$,  & $l\ge2$     & $R(\pi_2)$                   \\
  3 & $A_{2l}$,  & $l\ge2$     & $R(\pi_2)\oplus R(\pi_2)$    \\
  4 & $A_{2l}$,  & $l\ge2$     & $R(\pi_1)\oplus R(\pi_2)^*$  \\
  5 & $C_l$,     & $l\ge2$     & $R(\pi_1)$                   \\
  6 & $D_5$      &             & $R(\pi_4)$                   \\
\hline
\end{tabular}
\end{align*}
Here $\pi_k$ is a $k$-th fundamental weight
of  $\Lg$ (in the numbering of~\cite{VinO}).
$R(\lambda)$~is an irreducible representation
of $\Lg$ with the highest weight $\lambda$, $R'(\lambda)$~is its dual.
$R_1\oplus R_2$~is a direct sum of $R_1$ and~$R_2$,
$kR$~is a direct sum of $k$ copies of $R$.

The adjoint representation of a simple group has a unique
non-zero orbit of minimal dimension, namely,
the orbit of a highest root vector.
Therefore in the case~5 the only possible variant
was considered above.

It follows from~\eqref{hrep} that the centralizer
of $e$ in $\Lg$ coincides with the
isotropy subalgebra of the vector $\ol h$ 
with respect to the representation $\tG$. 
The centralizer of any non-zero element in a Lie algebra obviously has
a non-trivial center.
But it follows from~\cite{Elashvilli}
that the center of generic stabilizer is trivial
in the cases~2,~6, and 1 as $l>1$. If $l=1$ then the case~1 is 
contained in the case~5 and was considered above.

Consider the case~3.
By dimension count the nilpotent element
$e\in\Lg\iso\Lsl_n$ should be principal.
Therefore, the $\Lsl_2$-module $\Lg/\Lh$ is multiplicity--free.
But the  $\Lsl_2$-module $R(\pi_2)\oplus R(\pi_2)$
is obviously not multiplicity--free and we are done.

Consider the case~4.
By dimension count the nilpotent element
$e\in\Lg\iso\Lsl_n$ 
should be not principal.
Hence, the $\Lsl_2$-module
$\bbC^n$ is not irreducible.
Therefore, there exists a 1-PS $T\subset SL_n$
which commutes with $\Lsl_2$ (in particular, $T\subset H$) and 
does not have non-zero fixed vectors in $\bbC^n$. 
But $H$ is a centralizer of a generic element
in $R(\pi_1)\oplus R(\pi_2)^*$ and so 
$T$ should have fixed vectors
and we are done. The theorem is proved.

\end{document}